\documentclass{article}
\usepackage{graphicx,float,subfig,amsthm}               
\usepackage{amsmath,pgf,tikz,tkz-graph}
\usetikzlibrary{arrows}
\usepackage{graphicx}
\usepackage{amssymb}
\usepackage{epstopdf}
\usepackage[noblocks]{authblk}
\newtheorem{theorem}{Theorem}[section]
\newtheorem{lemma}[theorem]{Lemma}

\newtheorem{corollary}[theorem]{Corollary} 
\renewcommand{\qed}{\hfill $\blacksquare$}
\title{Constructing a Family of 4-Critical Planar Graphs with High Edge-Density}
\author[1]{Yao Tianxing}
\author[2]{Zhou Guofei\thanks{gfzhou@nju.edu.cn}}
\affil[1]{Sanjiang University, Nanjing, 210012, China}
\affil[2]{Department of Mathematics, Nanjing University, Nanjing 210093, China}
\date{}                                                       
\begin{document}
\baselineskip 0.23in
\maketitle
\begin{abstract}
A graph $G=(V,E)$ is a $k$-critical graph if $G$ is not $(k -1)$-colorable but  $G-e$ is $(k-1)$-colorable for every $e\in E(G)$. In this paper, we construct a family of 4-critical planar graphs with $n$ vertices and $\frac{7n-13}{3}$ edges. As a consequence, this improved the bound for the maximum edge density obtained by Abbott and Zhou.  We conjecture that this is the largest edge density for a 4-critical planar graph.
\end{abstract}
\section{Introduction}
Let $G=(V,E)$ be a graph, $G$ is said to be $k$-colorable if there is a assignment of $k$ colors to the vertices of $G$ such that no two adjacent vertices of $G$ get the same color.  The chromatic number of $G$, denoted by $\chi(G)$, is the least integer $k$ such that $G$ is $k$-colorable. 
A graph $G=(V,E)$ is a $k$-critical graph if $G$ is not $(k -1)$-colorable but  $G-e$ is $(k-1)$-colorable for every $e\in E(G)$. A planar graph is a graph that can be embedded in the plane, i.e., it can be drawn on the plane in such a way that its edges intersect only at their endpoints. Let $G=(V,E)$ be a graph and $v\in V(G)$, we denote by $N(v)$ the set of vertices that are adjacent to $v$. Let $f(n)$ and $F(n)$ denote respectively the minimum number and maximum number of a 4-critical planar graph with $n$ vertices. 
Table 1 shows some exact values of $f(n)$ and $F(n)$ for $n\le 14$.
In \cite{Yancey}, A.V. Kostochka and M. Yancey proved that $f(n)\ge \frac{5n-2}{3}$, and this bound is sharp in the sense that there are infinitely many 4-critical planar graphs on $n$ vertices and $\frac{5n-2}{3}$ edges. As for $F(n)$, H. L. Abbott and B. Zhou \cite{Abbott1} proved that $F(n)\le 2.75 n$, G. Koster \cite{GK} later  improved this bound to $\frac{5n}{2}$. We believe that this upper bound can not be obtained. Let $G=(V,E)$ be a graph, define $S={\rm sup} |E(G)|/|V(G)|$, where the bounds are taken over all 4-critical planar graphs with $|V(G)|$ vertices and $|E(G)|$ edges.
Gr\"{u}nbaum \cite{BG} used Haj\'os's construction \cite{GH} to show that $S\ge 79/39 = 2.02564...$; and he asked the question of determining the maximum edge density of planar 4-critical graphs. Abbott and Zhou \cite{Abbott1} used a variation of the Haj\'os's construction to show that $S\ge 39/19 = 2.05263...$. In this paper, we prove that $S\ge \frac{7}{3}$ by constructing a family of 4-critical planar graphs on $n$ vertices and $\frac{7n-13}{3}$ edges.

\begin{table}[h]
\renewcommand{\arraystretch}{1.3} 
\setlength{\tabcolsep}{0.4cm}
\centering
  \begin{tabular}{|c|c|c|c|c|c|c|c|c|c|}
  \hline
    $n$ & 6 & 7 & 8 & 9 & 10 & 11 & 12 & 13 & 14\\
    \hline
    $f(n)$ & 10 & 11 & 14 & 15 & 16 & 19 & 20 & 21 & 24\\
    \hline
    $F(n)$ & 10 & 12 & 14 & 16 & 18 & 20 & 22 & 26 & 28\\
    \hline
  \end{tabular}
  \caption{Some values of $f(n)$ and $F(n)$}
\end{table}

\section{Main Results}
\begin{figure}[H]
\centering
\begin{tikzpicture}[line cap=round,line join=round,>=triangle 45,x=1.0cm,y=1.0cm]
\clip(-4.3,0.2) rectangle (6.3,3.6);
\draw (-2.,3.)-- (0.,3.);
\draw (0.,3.)-- (2.,3.);
\draw (4.,3.)-- (6.,3.);
\draw (6.,3.)-- (5.,2.);
\draw (5.,2.)-- (4.,1.);
\draw (4.,1.)-- (3.,2.);
\draw (2.,3.)-- (1.,2.);
\draw (0.,1.)-- (1.,2.);
\draw (0.,1.)-- (-1.,2.);
\draw (-1.,2.)-- (-2.,3.);
\draw (-2.,3.)-- (-3.,2.);
\draw (-4.,1.)-- (-3.,2.);
\draw (-3.,2.)-- (-4.,3.);
\draw (-3.,2.)-- (-2.,1.);
\draw (-2.,1.)-- (-1.,2.);
\draw (-1.,2.)-- (0.,3.);
\draw (0.,3.)-- (1.,2.);
\draw (1.,2.)-- (2.,1.);
\draw (3.,2.)-- (4.,3.);
\draw (4.,3.)-- (5.,2.);
\draw (5.,2.)-- (6.,1.);
\draw (6.,1.)-- (4.,1.);
\draw (2.,1.)-- (0.,1.);
\draw (0.,1.)-- (-2.,1.);
\draw (-2.,1.)-- (-4.,1.);
\draw (-4.,3.)-- (-2.,3.);
\draw (-4.,3.5) node[anchor=north] {$u_0$};
\draw (-2.,3.5) node[anchor=north] {$u_1$};
\draw (0,3.5) node[anchor=north] {$u_2$};
\draw (2,3.5) node[anchor=north] {$u_3$};
\draw (4,3.5) node[anchor=north] {$u_{k-1}$};
\draw (6,3.5) node[anchor=north] {$u_{k}$};
\draw (-3,2.) node[anchor=east] {$w_1$};
\draw (-1,2.0) node[anchor=east] {$w_2$};
\draw (1.,2.) node[anchor=east] {$w_3$};
\draw (4.1,2.) node[anchor=east] {$w_{k-1}$};
\draw (5.7,2.) node[anchor=east] {$w_{k}$};
\draw (-4.,0.5) node[anchor=south] {$v_0$};
\draw (-2.,0.5) node[anchor=south] {$v_1$};
\draw (0.,0.5) node[anchor=south] {$v_2$};
\draw (2,0.5) node[anchor=south] {$v_3$};
\draw (4,0.5) node[anchor=south] {$v_{k-1}$};
\draw (6,0.5) node[anchor=south] {$v_{k}$};
\draw (2.6,3.) node[anchor=west] {$\cdots$};
\draw (1.8,2.0) node[anchor=west ] {$\cdots$};
\draw (2.6,1.) node[anchor=west] {$\cdots$};
\begin{scriptsize}
\draw [fill=black] (-4.,3.) circle (1.5pt);
\draw [fill=black] (-2.,3.) circle (1.5pt);
\draw [fill=black] (0.,3.) circle (1.5pt);
\draw [fill=black] (2.,3.) circle (1.5pt);
\draw [fill=black] (4.,3.) circle (1.5pt);
\draw [fill=black] (6.,3.) circle (1.5pt);
\draw [fill=black] (-3.,2.) circle (1.5pt);
\draw [fill=black] (-1.,2.) circle (1.5pt);
\draw [fill=black] (1.,2.) circle (1.5pt);
\draw [fill=black] (3.,2.) circle (1.5pt);
\draw [fill=black] (5.,2.) circle (1.5pt);
\draw [fill=black] (-4.,1.) circle (1.5pt);
\draw [fill=black] (-2.,1.) circle (1.5pt);
\draw [fill=black] (0.,1.) circle (1.5pt);
\draw [fill=black] (2.,1.) circle (1.5pt);
\draw [fill=black] (4.,1.) circle (1.5pt);
\draw [fill=black] (6.,1.) circle (1.5pt);
\end{scriptsize}
\end{tikzpicture}
\caption{The graph $H_k$}
\label{fig1}
\end{figure}
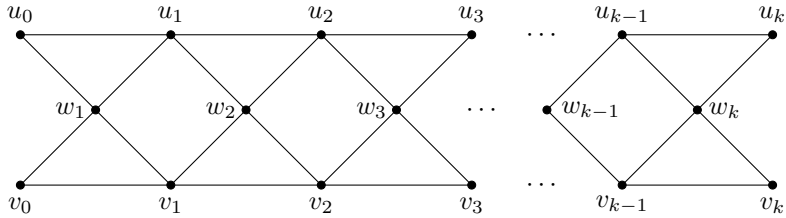
\begin{lemma}\label{lem1}
Let $H_k$ be the graph shown in Figure \ref{fig1}, then $H_k$ is three colorable; Moreover, let $c: V(H_k)\to \{1,2,3\}$ be a 3-coloring of $H_k$, 

(i) if there is nonnegative integer $i\ (0\le i\le k)$ such that $c(u_i)=c(v_i)$, then $c(u_j)=c(v_j)$ for all $j=0,1,2,\cdots,k$;

(ii) if there is nonnegative integer $i\ (0\le i\le k)$ such that $c(u_i)\ne c(v_i)$, then each $w_i$ gets the same third color, and $c(u_j)\ne c(v_j)$ for all $j=0,1,2,\cdots,k$;
\end{lemma}

{\bf Proof.} We color  $w_i$ ($1\le i\le k$) the color 1, and for each pair of vertices $u_{2i}$ and $v_{2i}$, we color them the color 2; finally, all the remaining vertices are colored 3. It is easy to see that it is a 3-coloring of $H_k$.

(i) It is easy to check that (i) is valid for $H_1$; Assume that  (i) holds for $H_{k-1}$. Now consider $H_k$, suppose there is nonnegative integer $i\ (0\le i\le k)$ such that $c(u_i)=c(v_i)$ (say $c(u_0)=c(v_0)$). Note that $H_{k-1}=H_k-\{u_k,v_k,w_k\}$, by assumption,  we have $c(u_j)=c(v_j)$ for all $j=0,1,2,\cdots,k-1$; 
without loss of generality, assume that $c(u_{k-1})=c(v_{k-1})=1$. Now, if $w_k=2$, then $c(u_{k-1})=c(v_{k-1})=3$; and if $w_k=3$, then $c(u_{k-1})=c(v_{k-1})=2$;  This proves (i) by induction.

(ii) We can prove it also by induction, the details are omitted here. \qed

\begin{theorem}\label{thm}
For each positive integer $k$, there is a 4-critical planar graph with  $6k+7$ vertices and $14k+12$ edges.
\end{theorem}

\begin{corollary}
$S\ge\frac{7}{3}$.
\end{corollary}

\begin{figure}[H]
\centering
\includegraphics[scale=1.1]{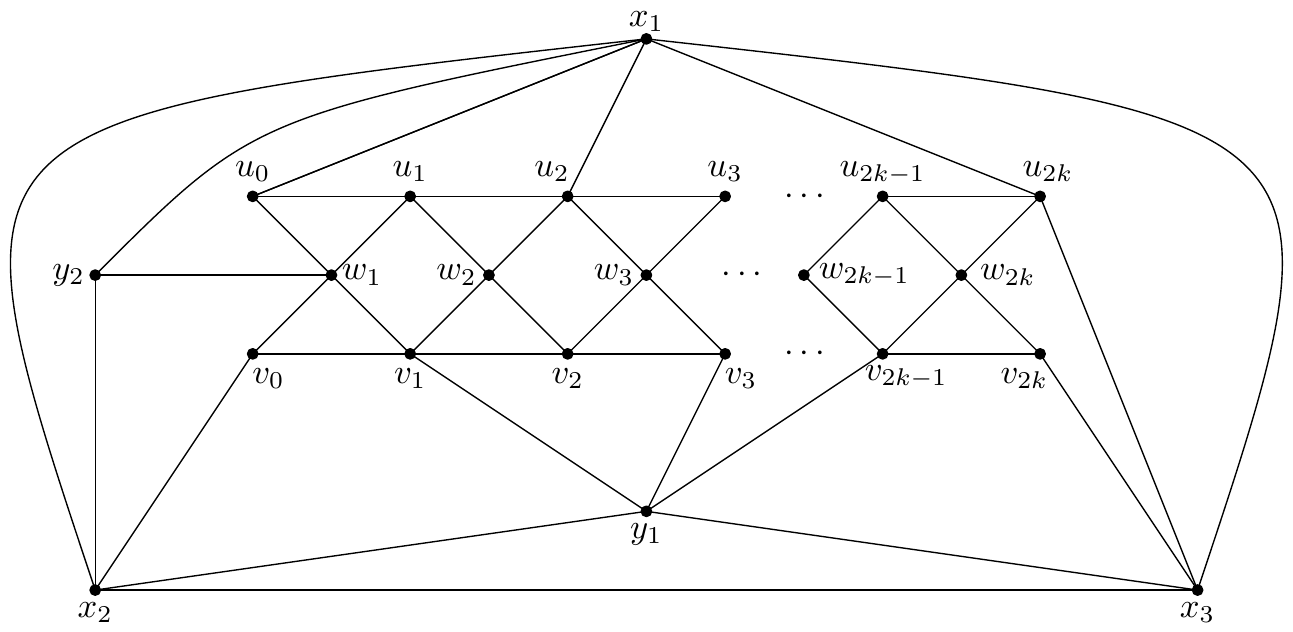}
\caption{The graph $G_k$}
\label{fig2}
\end{figure}

{\bf Proof of theorem \ref{thm}.} We construct a graph $G_k$ as shown in Figure \ref{fig2}, where $V(G_k)=V(H_{2k})\cup\{x_1,x_2,x_3,y_1,y_2\}$ and 
\begin{eqnarray*}
E(G_k)=E(H_{2k}) &\cup & \{x_1u_0,x_1u_2,\cdots,x_1u_{2k}\}\\
&\cup & \{y_1v_1,y_1v_3,\cdots,y_1v_{2k-1}\}\\
&\cup & \{x_1x_2,x_1x_3,x_1y_2,x_2y_1,x_2y_2,x_2v_{0},x_2x_3\}\\
&\cup & \{x_3y_1,x_3v_{2k},x_3u_{2k},y_2w_{1}\}
\end{eqnarray*}

\begin{figure}[H]
\centering
\subfloat[$G_1$]{
\includegraphics[width=0.4\textwidth]{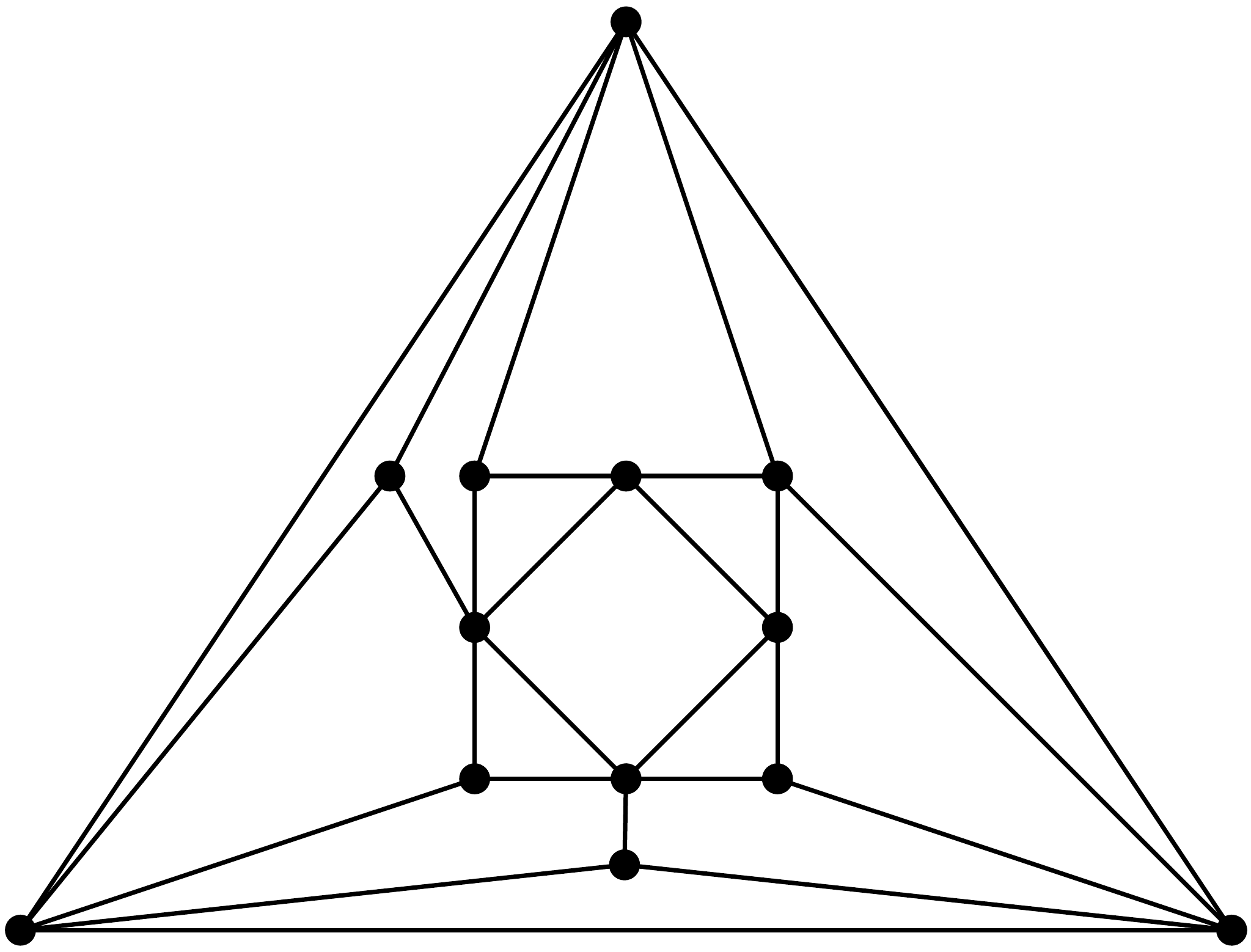}}
\hskip 2cm
\subfloat[$G_2$]{
\includegraphics[width=0.4\textwidth]{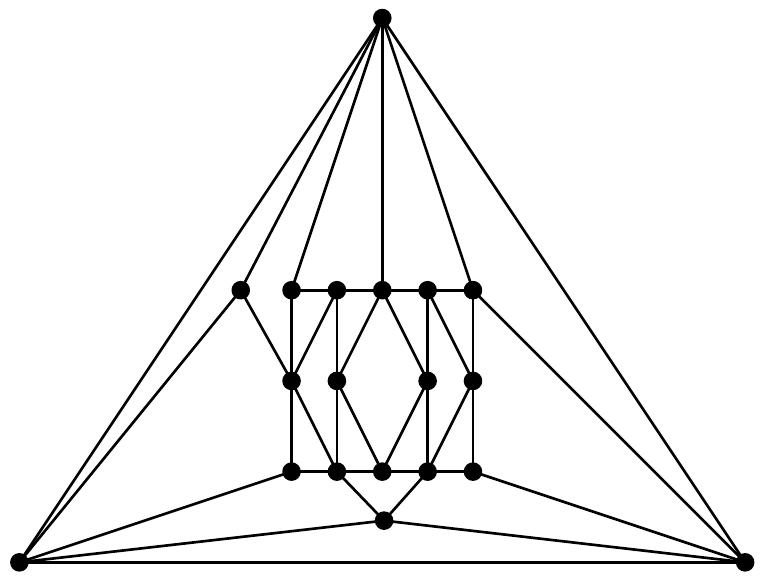}}
\caption{The 4-critical graphs $G_1$ and $G_2$}
\label{fig4}
\end{figure}

When $k=1$ and $k=2$, see Figure \ref{fig4} in particular, we can check that both of them are 4-critical planar graphs.

Note that $G_k$ is a planar graph with $6k+7$ vertices and $14k+12$ edges. In the following, we shall prove that 
both $G_k$ is 4-critical. 

First, we prove that $G_k$ is not 3-colorable.  Suppose that $G_k$ is 3-colorable, let $c: V(G_k)\to \{1,2,3\}$ be a 3-coloring of $G_k$. Note that $\{x_1,x_2,x_3\}$ form a triangle in $G_k$, without loss of generality, assume that $c(x_1)=1,c(x_2)=2,c(x_3)=3$, then $c(y_1)=1, c(y_2)=3, c(u_{2k})=2$. Note that $c(v_{2k})\in\{1,2\}$ since $x_3$ and $v_{2k}$ are adjacent. If $c(v_{2k})=c(u_{2k})=2$, then $c(v_{2k-1})=3$ since $v_{2k-1}$ is adjacent to both $y_1$ and $v_{2k}$. By Lemma \ref{lem1}, we have that $c(u_{2k-1})=c(v_{2k-1})=3$. By this way, we get in general that for each $u_i\in N(x_1)$, $c(u_i)=c(v_i)=2$; and for each $v_j\in N(y_1)$, $c(v_j)=c(u_j)=3$. So $c(v_{0})=2$, this is impossible since $c(x_2)=2$. If $c(v_{2k})\ne c(u_{2k})$, then $c(v_{2k})=1$. By Lemma \ref{lem1}, $c(w_i)=3$ for $1\le i\le 2k$. So $c(y_2)=c(w_{1})=3$, this is impossible since $y_2$ and $w_{1}$ are adjacent. Therefore, $G_k$ is not 3-colorable.

Next, we prove that for each $e\in E(G_k)$, $G_k-e$ is 3-colorable. Since the number of automorphisms of $G_k$ is one, $G_k$ is extremely not symmetric, so we have to consider awkwardly every possible edge of  $G_k$. In the following, let 

\hskip 3cm $U_1=N(x_1)\cap V(H_{2k});$

\hskip 3cm $U_2=\{u_0,u_1,u_2,\cdots,u_{2k}\}-U_1;$

\hskip 3cm $W=\{w_1,w_2,\cdots,w_{2k}\};$

\hskip 3cm $V_1=N(y_1)\cap V(H_{2k});$

\hskip 3cm $V_2=\{v_0,v_1,v_2,\cdots,v_{2k}\}-V_1;$

\hskip 3cm $[u_i,u_j] = \{u_i,u_{i+1},u_{i+2},\cdots,u_j\}$, where $i<j$;

\hskip 3cm $[u_i:u_j] = \{u_t\ |\ u_t\in[u_i,u_j] \ {\rm and}\ t-i\ {\rm is \ even}\}$, where $i<j$.

Furthermore, if a graph $G$ is 3-colorable, we will denote by $C_1$ and $C_2$ the set of vertices that are colored 1 and 2 respectively. For the sake of conciseness, we will not indicate the set of vertices that are colored 3.

Note that $G_k$ has four vertices $y_2,u_0,v_0,v_{2k}$ which have degree 3. If $G_k-v$ is 3-colorable for $v\in\{y_2,u_0,v_0,v_{2k}\}$, then it is obvious that $G_k-e$ is 3-colorable for each edge $e$ that is incident with $v$. Let $H=G_k-v$ for some $v\in \{y_2,u_0,v_0,v_{2k}\}$, we first prove that $H$ is 3-colorable.

If $H=G_k-y_2$, let $C_1=\{x_1,y_1\}\cup U_2\cup V_2$, $C_2=\{x_3\}\cup W$;

If $H=G_k-u_0$, let $C_1=\{x_1,y_1,v_0\}\cup [w_2,w_{2k}]$, $C_2=\{x_3,y_2\}\cup U_2\cup V_1$;

If $H=G_k-v_0$, let $C_1=\{x_1,y_1\}\cup W$, $C_2=\{x_3,y_2\}\cup U_2\cup V_1$;

If $H=G_k-v_{2k}$, let $C_1=\{x_1,y_1\}\cup W$, $C_2=\{x_2\}\cup U_1\cup V_1$.

This proves that $H=G_k-v$ is 3-colorable for some $v\in \{y_2,u_0,v_0,v_{2k}\}$.

Next, let $e$ be an edge that is not incident with any vertex of $\{y_2,u_0,v_0,v_{2k}\}$ and let $G=G_k-e$, we shall prove that $G$ is 3-colorable.

If $e=x_1x_2$, let $C_1=\{x_1,x_2\}\cup W$, $C_2=\{y_1\}\cup U_1\cup V_2$;

If $e=x_1x_3$, let $C_1=\{x_1,x_3\}\cup W$, $C_2=\{y_1,y_2\}\cup U_1\cup V_2$;

If $e=x_2x_3$, let $C_1=\{x_1\}\cup W$, $C_2=\{y_1,y_2\}\cup U_1\cup V_2$;

If $e=x_2y_1$, let $C_1=\{x_1\}\cup U_2\cup V_1$, $C_2=\{x_3,y_2,u_0,v_0\}\cup [w_2,w_{2k}]$;

If $e=x_3y_1$, let $C_1=\{x_1\}\cup U_2\cup V_1$, $C_2=\{x_3,y_1,y_2,u_0,v_0\}\cup [w_2,w_{2k}]$;

If $e=x_3u_{2k}$, let $C_1=\{x_1,y_1\}\cup U_2\cup V_2$, $C_2=\{x_2\}\cup W$;

If $e=x_1u_{i}$ ($i$ is even), let $C_1=\{x_1,y_1,u_i,v_i\}\cup W-\{w_i,w_{i+1}\}$, $C_2=\{x_3,y_2,w_i\}\cup [u_0:u_{i-2}]\cup [u_{i+1}:u_{2k-1}]\cup [v_0:v_{i-2}]\cup [v_{i+1}:v_{2k-1}]$;

If  $e=y_1v_{i}$ ($i$ is odd), let $C_1=\{x_1,y_1,u_i,v_i\}\cup W-\{w_i,w_{i+1}\}$, $C_2=\{x_2,w_i\}\cup [u_1:u_{i-2}]\cup [u_{i+1}:u_{2k}]\cup [v_1:v_{i-2}]\cup [v_{i+1}:v_{2k}]$;

If  $e=u_iu_{i+1}$ ($i$ is odd), let $C_1=\{x_1,y_1\}\cup [u_{i+2}:u_{2k-1}]\cup [v_{i+1}:v_{2k}]\cup[w_1,w_i]$, $C_2=\{x_2\}\cup [u_1:u_{i}]\cup [u_{i+1}:u_{2k}]\cup V_1$;

If  $e=u_iu_{i+1}$ ($i$ is even), let $C_1=\{x_1,y_1\}\cup [u_1:u_{i-1}]\cup [v_{0}:v_{i}]\cup[w_{i+2},w_{2k}]$, 
$C_2=\{x_3,y_2\}\cup [u_0:u_{i}]\cup [u_{i+1}:u_{2k-1}]\cup V_1$;

If  $e=v_iv_{i+1}$ ($i$ is even), let $C_1=\{x_1,y_1\}\cup W$, 
$C_2=\{x_3,y_2\}\cup [v_0:v_{i}]\cup [v_{i+1}:v_{2k-1}]\cup U_2$;

If  $e=v_iv_{i+1}$ ($i$ is odd), let $C_1=\{x_1,y_1\}\cup W$, 
$C_2=\{x_3,y_2\}\cup [v_0:v_{i-1}]\cup [v_{i+2}:v_{2k-1}]\cup U_2$;

If  $e=u_iw_{i}$ ($u_i\in U_2$), let $C_1=\{x_1,y_1\}\cup [u_i:u_{2k-1}]\cup [v_{i+1}:v_{2k}]\cup [w_1,w_i]$, 
$C_2=\{x_2\}\cup U_1\cup V_1$;

If  $e=u_iw_{i+1}$ ($u_i\in U_2$), let $C_1=\{x_1,y_1\}\cup [u_{i+2}:u_{2k-1}]\cup [v_{i+1}:v_{2k}]\cup [w_1,w_i]$, 
$C_2=\{x_2\}\cup U_1\cup V_1$;

If  $e=u_iw_{i+1}$ ($u_i\in U_1$), let $C_1=\{x_1,y_1\}\cup U_2\cup V_2$, 
$C_2=\{x_3,y_2\}\cup [u_0:u_i]\cup[v_1:v_{i-1}]\cup [w_{i+1},w_{2k}]$;

If  $e=u_iw_{i}$ ($u_i\in U_1$), let $C_1=\{x_1,y_1\}\cup U_2\cup V_2$, 
$C_2=\{x_3,y_2\}\cup [u_0:u_{i-1}]\cup[v_1:v_{i-1}]\cup [w_{i+1},w_{2k}]$;

If  $e=v_{i+1}w_{i+1}$ ($v_{i+1}\in V_1$), let $C_1=\{x_1,y_1\}\cup U_2\cup V_2$, 
$C_2=\{x_3,y_2\}\cup [u_0:u_{i-1}]\cup[v_1:v_{i-1}]\cup [w_{i+3},w_{2k}]$;

If  $e=v_{i-1}w_{i}$ ($v_{i-1}\in V_1$), let $C_1=\{x_1,y_1\}\cup U_2\cup V_2$, 
$C_2=\{x_3,y_2\}\cup [u_0:u_{i-1}]\cup[v_1:v_{i-1}]\cup [w_{i},w_{2k}]$;

If  $e=v_{i}w_{i}$ ($v_{i}\in V_2$), let $C_1=\{x_1,y_1\}\cup [u_{i+1}:u_{2k-1}]\cup [v_i: v_{2k}]\cup [w_1,w_i]$, 
$C_2=\{x_2\}\cup U_1\cup V_1$;

If  $e=v_{i}w_{i+1}$ ($v_{i}\in V_2$), let $C_1=\{x_1,y_1\}\cup [u_{i+1}:u_{2k-1}]\cup [v_{i+2}: v_{2k}]\cup [w_1,w_i]$, 
$C_2=\{x_2\}\cup U_1\cup V_1$.

From the above coloring schedules, we have checked that for each $e\in E(G_k)$, $G_k-e$ is 3-colorable. This completes the proof of theorem \ref{thm}. \qed

\section{Some remarks and problems}
\begin{itemize}
\item Note that both $G_k$ and $G_k^2$ have minimum degree 3. Is there a 4-critical planar graph on  $6k+7$ vertices and $14k+12$ edges and with $\delta \ge 4$?
\item Note that in Table 1, all the values of $F(n)$ for small $n$ are even. Is it true that $F(n)$ is even for all positive integer $n$?
\item We conjecture that $S=\frac{7}{3}$. More concisely, we conjecture that $F(n)\le \frac{7n-13}{3}$, where the equality holds only if $n\equiv 1({\rm mod}\ 6)$.
\end{itemize}

\noindent{\bf Acknowledgments}

This research work is supported by National Natural Science Foundation of China under grant No. 11571168 and 11371193.

\end{document}